\documentclass[11pt]{article}
\usepackage{amsmath,verbatim,amssymb,array}
\usepackage[cp1250]{inputenc}
\usepackage[T1]{fontenc}
\usepackage[english]{babel}
\usepackage{a4wide,cite}
\usepackage{graphicx}
\usepackage{bm}
\usepackage{color}

\newenvironment{pf}{{\it Proof.\ }\rm }{ \hfill$\Box$\ \\[0.25ex]}
\numberwithin{equation}{section}
\newtheorem{thm}{Theorem}[section]
\newtheorem{lem}[thm]{Lemma}
\newtheorem{cor}[thm]{Corollary}
\newtheorem{alg}[thm]{Algorithm}
\newtheorem{prob}[thm]{Problem}

\newenvironment{rem}{\refstepcounter{thm}\ \\[2mm]%
                  \noindent{\it Remark~\thethm\ \ }}{\ \\[2mm]\rm}

\newcommand{\intT}{\mbox{$\displaystyle\int\!\!\!\!\int\limits_{\!\!\!\!T}$}}
\newcommand{\N}{\mathbb N}

\newcommand{\R}{\mathbb R}
\newcommand{\Z}{\mathbb Z}
\newcommand{\Zplus}{\Z_+}
\newcommand{\C}{\mathbb C}

\newcommand{\sbm}[1]{\mbox{\scriptsize$\bm{#1}$}}  
\newcommand{\tbm}[1]{\mbox{\tiny$\bm{#1}$}} 
\newcommand{\bma}{\mbox{$\bm{\alpha}$}}
\newcommand{\lbma}{\mbox{$\left|\bma\right|$}}
\newcommand{\sbma}{\mbox{\scriptsize$\bma$}}

\newcommand{\tbma}{\mbox{$\tilde{\bma}$}}
\newcommand{\tsbm}[1]{\mbox{\scriptsize$\tilde{\bm{#1}}$}}
\newcommand{\tsbma}{\mbox{\scriptsize$\tsbm{\alpha}$}}

\newcommand{\bmc}{\mbox{$\bm c$}}
\newcommand{\sbmc}{\mbox{$\sbm c$}}
\newcommand{\lbmc}{\mbox{$\left|\bmc\right|$}}

\newcommand{\dx}{\,\mbox{{\rm d}$\bm x$}}

\newcommand{\sP}{\mbox{$\mathsf P$}}
\newcommand{\sR}{\mbox{$\mathsf R$}}
\newcommand{\sQ}{\mbox{$\mathsf Q$}}

\newcommand{\sD}{\mbox{$\mathcal D$}}
\newcommand{\cR}{\mbox{$\mathcal R$}}

\newcommand{\rbinom}[2]{\mbox{$\displaystyle\binom{#1}{#2}^{\!\!-1}\!\!$}}
\newcommand{\hyper}[5]{\,{}_{#1}F_{#2}\!\left(%
				\begin{array}{cc}{#3}\\[0.35ex]{#4} \end{array}\Big|\,{#5} \right)}

\begin{document}
	\title{\LARGE\bf B\'ezier form of dual bivariate Bernstein polynomials}
\author{\normalsize \textbf{Stanis{\l}aw Lewanowicz}${}^{a,\!}$%
           	 \footnote{Corresponding author\newline
           	 	 \hspace*{1.5em}\textit{Email addresses}: Stanislaw.Lewanowicz@cs.uni.wroc.pl 
		(Stanisław Lewanowicz),\newline
		 Pawel.Keller@mini.pw.edu.pl (Paweł Keller), Pawel.Wozny@cs.uni.wroc.pl (Paweł Woźny)}
		\,\:, \textbf{Pawe{\l} Keller}${}^{b}$, \textbf{Pawe{\l} Wo\'zny}${}^{a}$}

\date{{\small\it ${}^{a}\!\!\!$ Institute of Computer Science, University of Wroc{\l}aw,
 	 ul.~Joliot-Curie 15, 50-383 Wroc{\l}aw, Poland\\[0.5ex]	        
 	 \noindent ${}^{b}\!\!\!$ Faculty of Mathematics and Information Science, 
 	Warsaw University of Technology,\\ ul. Koszykowa 75, 00-662 Warszawa, Poland}}
   
\maketitle		
 \thispagestyle{empty}  
 		 
\noindent {\small\textbf{Abstract}
   Dual Bernstein polynomials of one or two variables have proved to be very useful
	in obtaining B\'{e}zier form of the $L^2$-solution of the problem  of best polynomial
	approximation of B\'{e}zier curve or surface. 
	In this connection, the  B\'{e}zier coefficients of dual 
	Bernstein polynomials are to be evaluated at a reasonable cost.
	In this paper, a set of recurrence relations satisfied by the B\'{e}zier coefficients of dual 
	bivariate Bernstein polynomials is derived and an efficient algorithm 
	for evaluation of these coefficients is proposed.
        Applications of this result to some approximation problems of Computer Aided Geometric Design
(CAGD) are discussed.}\\[0.5ex]	 

\section{Introduction and preliminaries}                        
					\label{S:intro} 

The Bernstein-B\'ezier curves and surfaces have become the standard in the CAGD context
due to their favourable geometric properties (cf.\ \cite{Far02}).
Degree reduction, which consists in approximating a B\'ezier curve or surface by another one
of a lower degree, and  approximation of a  rational B\'ezier
curve or surface by a polynomial one have many important applications in geometric
modelling, such as data exchange, data compression, and data comparison.
There have been many papers relevant to this problem (see, e.g., 
\cite{Ahn03,ALPY04,CW02,Eck95,Hu13,LBPY02,LuW06,LuW08,%
		Rab05,RLY06,RA06,Sun05,SL04,Sha13,WL09,WL10}). 

The above approximation problems, with or without constraints,
are studied for different choices of the error norm.
Notice that the Bernstein polynomials do not form an orthogonal base,
so obtaining  the B\'ezier form of  the best  $L_2$-norm solution  
 requires some effort. In a frequently used approach, 
the main tool applied was transformation between the Bernstein and
orthogonal polynomial bases. Such methods are not only expensive, but also may be ill-conditioned
(cf.\ a remark in \cite{LuW08}).

Recently \cite{WL09,WL10}, a novel approach to the best $L_2$ approximation problem, 
using the dual basis associated with the univariate or bivariate Bernstein basis,
was proposed. The new methods do not use basis transformation matrices explicitly,
hence they do not share the abovementioned limitation.
High efficiency of the methods was obtained thanks to the application of recursive properties 
of the dual Bernstein polynomials.  
 
In the present paper, we derive the recurrence relations satisfied by the B\'{e}zier coefficients
of the dual bivariate Bernstein polynomials
and prove in detail a low-cost algorithm for numerical evaluation of these coefficients; 
see Section~\ref{S:EE}.
This algorithm (without proof) was successfully used in \cite{LKW17} as a part of the method for polynomial
approximation of rational triangular B\'{e}zier surfaces%
\footnote{Notice that for a specific choice of parameters, the rational B\'{e}zier surface reduces 
	to  a polynomial B\'{e}zier surface, so that the above approximation problem, 
with an additional assumption, is actually  the polynomial degree reduction problem discussed 
in \cite{WL10}. In the method proposed there, dual bivariate Bernstein polynomials also play 
the basic role. It should be stressed, however, that the present approach is substantially different
and much more efficient than the one proposed in \cite{WL10}.}; see the discussion in
Section~\ref{S:Appl}. Several useful properties of the dual bivariate Bernstein basis polynomials,
exploited in the proposed method, are obtained in Appendix~A. 
In Appendix~B, some results on the  bivariate Jacobi and Hahn  orthogonal polynomials are collected. 

Below, we introduce  notation and definitions that are used in the paper.

\begin{figure}[htb]
	\begin{center}     
        \begin{tabular}{ccc}
	\hspace*{-1cm} 
	\includegraphics[scale=0.225,angle=-90]{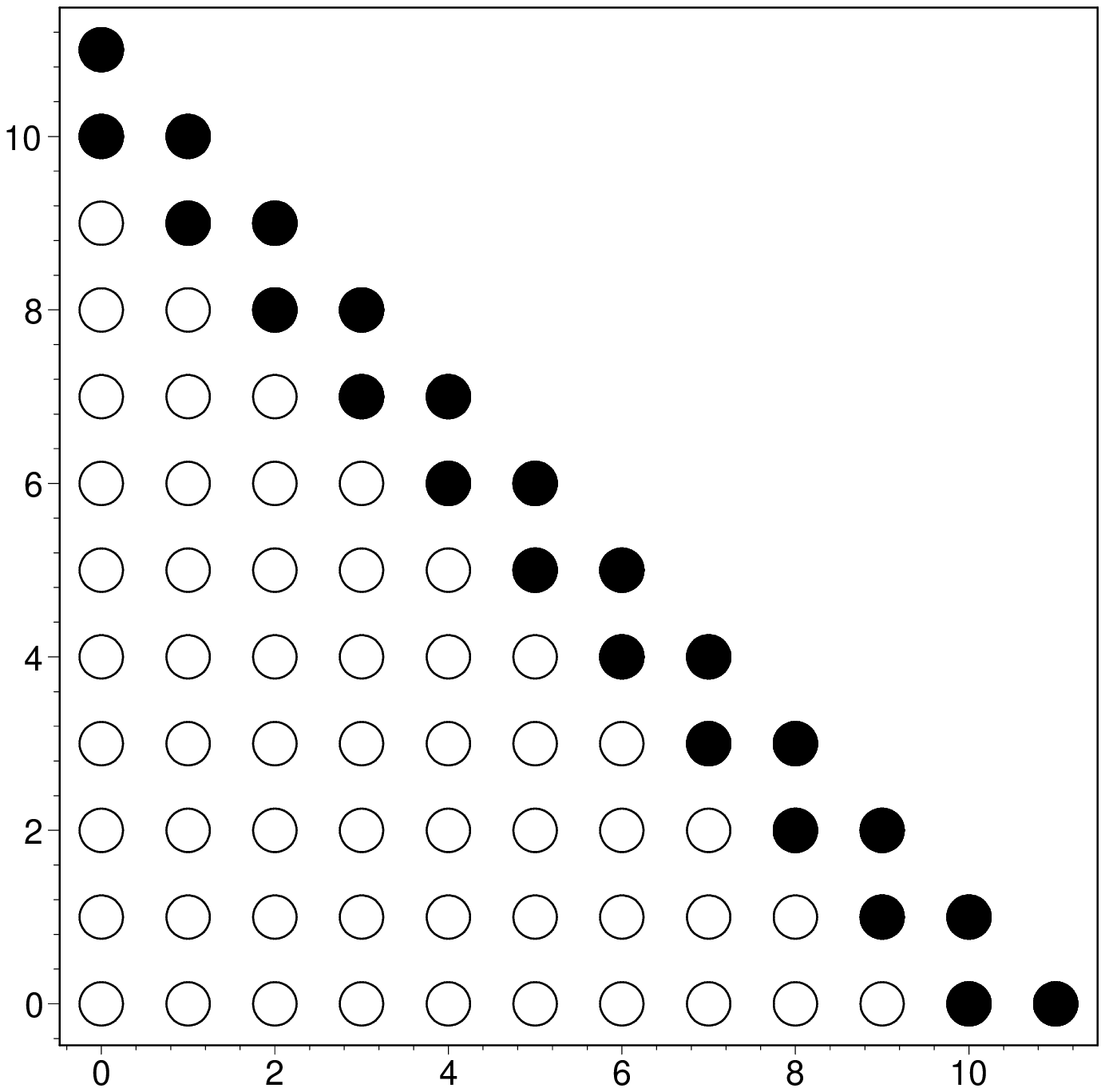} &
	\hspace*{-1.8cm} 
	\includegraphics[scale=0.225,angle=-90]{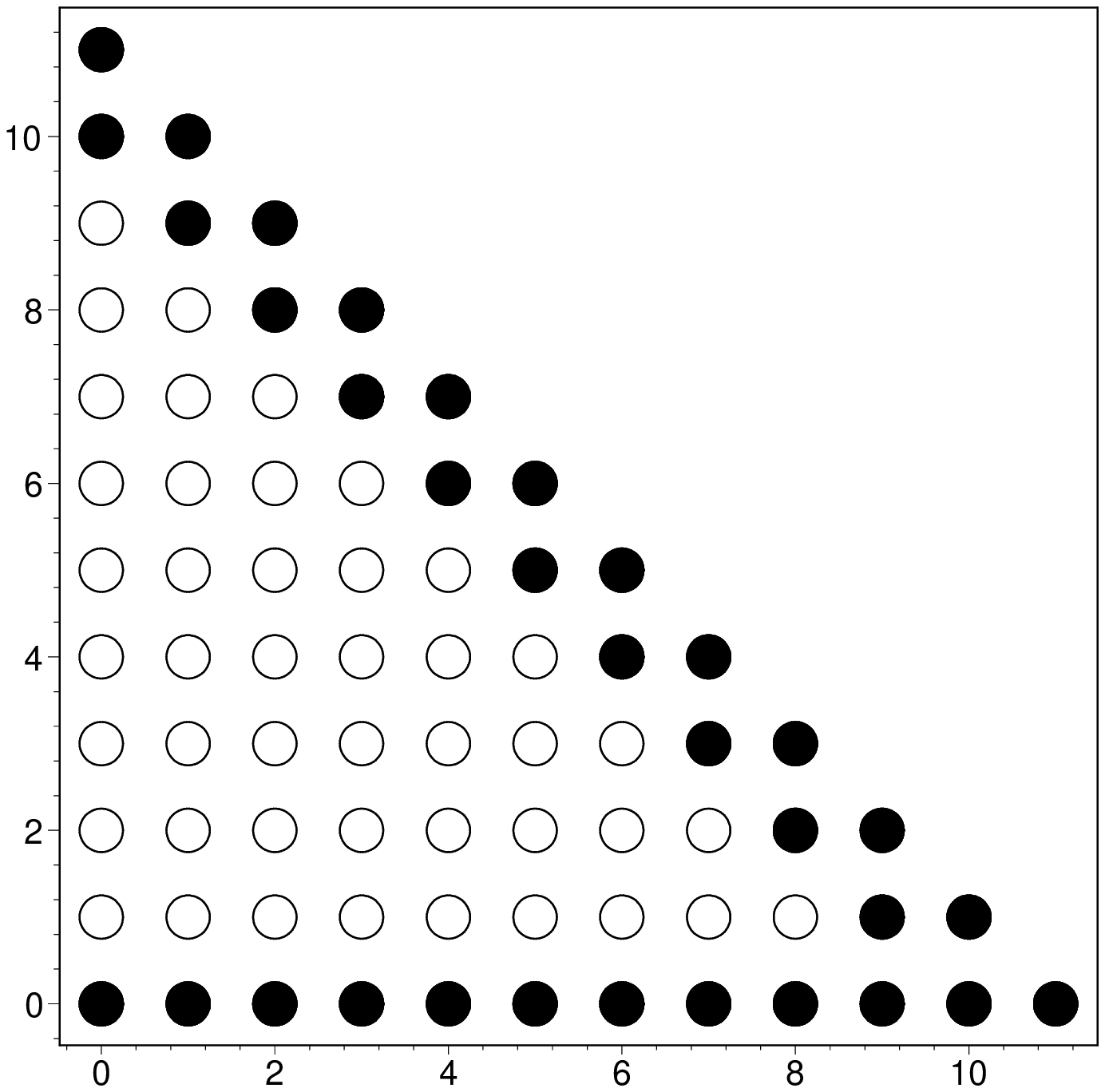}&
	\hspace*{-1.8cm} 
	\includegraphics[scale=0.225,angle=-90]{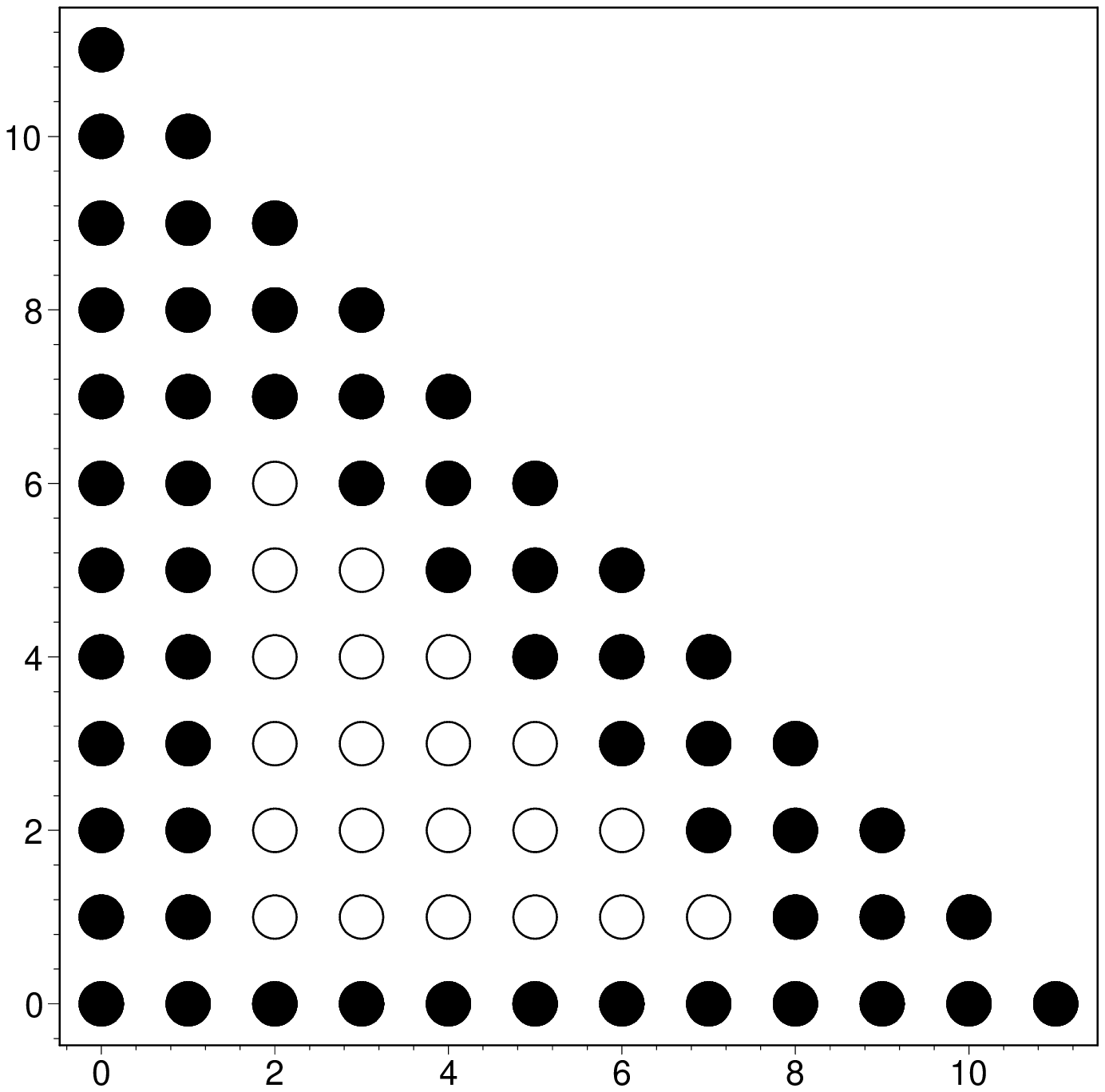}
	\\[-0.2ex]
	\hspace*{-0.9cm} $\bm{c}=(0,0,2)$ & \hspace*{-1.65cm} $\bm{c}=(0,1,2)$ 
	& \hspace*{-1.65cm}	$\bm{c}=(2,1,3)$ 
	\end{tabular}	 		
	\vspace*{0.125cm}
	\caption{Examples of sets \eqref{E:TOG} ($n=11$). Points of the set
		$\Omega^{\sbm c}_n$ 
	are marked by white discs, while the points of the set
	$\Gamma^{\sbm c}_n$ -- by black discs. Obviously, 
	$\Theta_n=\Omega^{\sbm c}_n\cup\Gamma^{\sbm c}_n$.%
	\label{fig:Fig-bc}}
	\end{center}
	\end{figure}

For  $\bm y:=(y_1,y_2,\ldots,y_d)\in\R^d$, we  denote 
\[
	|\bm y|:=y_1+y_2+\ldots+y_d,\qquad  \|\bm y\|:=\left(y^2_1+y^2_2+\ldots+y^2_d\right)^{\frac12}.
\]	
For $n\in\N$ and  $\bm c:=(c_1,\,c_2,\,c_3)\in\N^3$ such that $\lbmc<n$, we define 
the following sets (cf.\ Figure~\ref{fig:Fig-bc}):
\begin{equation}\label{E:TOG}
	\left.\begin{array}{l}
	\Theta_n:=\{\bm k=(k_1,k_2)\in\N^2:\: 0\le|\bm k|\le n\},  \\[1ex]
	\Omega^{\sbm c}_n:=\{\bm k=(k_1,k_2)\in\N^2:\:k_1\ge c_1,\,k_2\ge c_2,\,|\bm k|\le n-c_3\},\\[1ex]
	\Gamma^{\sbm c}_n:=\Theta_n\setminus\Omega^{\sbm c}_n.
	  \end{array}\;\right\}
\end{equation}	

\begin{rem}\label{R:TG}
		In applications to approximation problems related to  triangular
		(rational or polynomial) B\'ezier patches, 
		the set $\Theta_n$ corresponds to the set of control points of a  patch, while its subset $\Gamma_n^{\sbmc}$ \--- to the boundary points, where some constraints are 
		to be imposed.
		See Section~\ref{S:Appl}.
\end{rem}

Let $T$ be the standard triangle in $\R^2$,
\begin{equation}
	\label{E:T}
		T:=\{(x_1,\,x_2)\,:\,x_1,\,x_2\ge0,\: x_1+x_2\le1\}.		
\end{equation}	
For $n\in\N$ and $\bm k:=(k_1,k_2)\in\Theta_n$, we denote
\[
	\binom{n}{\bm k}:=\frac{n!}{k_1!k_2!(n-|\bm k|)!}.
\]
The \textit{shifted factorial} is defined for any $a\in\C$ by
\[
	(a)_0:=1;\quad (a)_k:=a(a+1)\cdots(a+k-1),\quad   k\ge1.
\]
The \textit{Bernstein polynomial basis} in $\Pi^2_n$, $n\in\N$, is given by
(see, e.g., \cite{Far86}, or \cite[\S 18.4]{Far02})
\begin{equation}\label{E:Ber2}
	B^n_{\sbm k}(\bm x):=\binom{n}{\bm k}x_1^{k_1}x_2^{k_2}(1-|\bm x|)^{n-|\sbm k|},
	\qquad \bm k:=(k_1,k_2)\in\Theta_n,\quad\bm x:=(x_1,x_2).
\end{equation}

The (unconstrained) \textit{dual bivariate Bernstein basis polynomials} \cite{LW06}, 
\begin{equation}\label{E:dBer2}
	D^n_{\sbm k}(\bm\cdot;\bma)\in\Pi^2_n,	
	\qquad \bm k\in\Theta_n,
\end{equation}
are defined so that
\[  
	\left\langle D^{n}_{\sbm k},\,B^n_{\sbm l} \right\rangle_{\sbma}=\delta_{\sbm k, \sbm l},
	\qquad \bm k,\bm l\in\Theta_n.
\] 
Here, $\delta_{\sbm k, \sbm l}$ equals 1 if $\bm k= \bm l$,  and 0 otherwise, while the inner product is defined by
\begin{equation}\label{E:Jinprod}
\langle f, g \rangle_{\sbma} :=\intT w_{\sbma}(\bm x)f(\bm x)\,g(\bm x)\dx,
\end{equation} 
where the weight function $w_{\sbma}$ is given by
\begin{equation}
	\label{E:w}
	w_{\sbma}(\bm x):=A_{\sbma}x_1^{\alpha_1}x_2^{\alpha_2}(1-|\bm x|)^{\alpha_3}, 
	\qquad\bma:=(\alpha_1,\,\alpha_2,\,\alpha_3), \quad\alpha_i>-1,
		\end{equation}	
and $A_{\sbma}$ is a normalisation factor (see \eqref{E:A}). 

In the sequel,  we use the notation $e^{\sbm k}_{\sbm l}(\bma,n)$ for the \textit{connection coefficients}
between bivariate Bernstein and  dual bivariate  Bernstein bases such that 
\begin{equation}\label{E:D2inB2-in}
	D^n_{\sbm k}(\bm x;\bma)=\,
	\sum_{\sbm l\in\Theta_n}\,e^{\sbm k}_{\sbm l}(\bma,n)\,
	B^n_{\sbm l}(\bm x),\qquad \bm k\in\Theta_n.
\end{equation}
Investigating the properties of these coefficients is one of the main goals of the paper.
The fast recursive scheme for computing the coefficients $e^{\sbm k}_{\sbm l}(\bma,n)$,
$\bm k,\bm l\in\Theta_n$, is formulated in the next section, while the wide
mathematical background is given in Section~\ref{SS:dualB2}.

For $n\in\N$ and  $\bm c:=(c_1,\,c_2,\,c_3)\in\N^3$ such that $\lbmc<n$, define 
 the constrained bivariate  polynomial space 
\[  
	\Pi^{\sbm c,\,2}_n :=\left\{P\in\Pi^2_n\::\: 
	P(\bm x)= x_1^{c_1}x_2^{c_2}(1-|\bm x|)^{c_3}\cdot\,Q(\bm x),
\;\mbox{where}\;Q\in\Pi^2_{n-|\sbm c|}\right\}.
\] 
It can be easily seen that the constrained set 
of  polynomials $B^n_{\sbm k}$ with the range of $\bm k$ restricted to 
$\Omega^{\sbm c}_n$ forms a basis in this space. 
We define the \textit{constrained  dual bivariate Bernstein basis polynomials}, 
\begin{equation}\label{E:constrdBer2}
	D^{(n,\sbm c)}_{\sbm k}(\bm{\cdot};\bma)\in\Pi^{\sbm c,\,2}_n,
	\qquad \bm k\in\Omega^{\sbm c}_n,
\end{equation}
 so that
\[ 
	\left\langle D^{(n,\sbm c)}_{\sbm k},\,B^n_{\sbm l}
	\right\rangle_{\sbma}=\delta_{\sbm k, \sbm l} \quad \mbox{for}
	\quad \bm k,\bm l\in\Omega^{\sbmc}_n,
\]
where the notation of \eqref{E:Jinprod} is used. For $\bm c=(0,0,0)$, basis \eqref{E:constrdBer2}
reduces to the unconstrained basis \eqref{E:dBer2} in $\Pi^2_n$. Obviously, for any $\sQ\in\Pi^{\sbm c,\,2}_n$,
we have
\begin{equation}\label{E:Bezform}
	\sQ(\bm x)=\sum_{\sbm k\in\Omega^{\sbmc}_n}\left\langle \sQ, D^{(n,\sbm c)}_{\sbm k}\right\rangle_{\sbma}\,B^n_{\sbm k}(\bm x).
\end{equation}	

Let $E^{\sbm k}_{\sbm l}(\bma,\bmc,n)$, $\bm l\in\Omega^{\sbm c}_n$, denote the B\'ezier  coefficients 
of the constrained bivariate dual Bernstein polynomial $D^{(n,\sbm c)}_{\sbm k}(\bm x;\bma)$,
	$\bm k\in\Omega^{\sbm c}_n$:
\begin{equation}
	\label{E:Dc2inB2-in}
	D^{(n,\sbm c)}_{\sbm k}(\bm x;\bma)=\sum_{\sbm l\in\Omega^{\sbm c}_n}
	\,E^{\sbm k}_{\sbm l}(\bma,\bmc,n)\,B^n_{\sbm l}(\bm x).
\end{equation}
According to Lemma \ref{L:Dc2inB2}, the following formula  relates the above coefficients
with the ''unconstrained'' connection coefficients defined in \eqref{E:D2inB2-in}:
\begin{equation}
	\label{E:Dc2inB2-coef-impl}
	E^{\sbm k}_{\sbm l}(\bma,\bmc,n):=U\,V_{\sbm k}(n)\,V_{\sbm l}(n)\,
	\,e^{\sbm k-\sbm c'}_{\sbm l-\sbm c'}(\bma+2\bm c,n-|\bm c|),\qquad {\bm k}, {\bm l}\in\Omega^{\sbm c}_n,
	\end{equation}
where $\bm c':=(c_1,c_2)$, and $U,V_{\sbm k}(n)$, $\bm k\in\Omega^{\sbm c}_n$, 
are constants defined in \eqref{E:UV}.

Hence, for fast evaluation of the coefficients \eqref{E:Dc2inB2-coef-impl}, it would be sufficient
to have  an efficient algorithm for computing
the  quantities $e^{\sbm k}_{\sbm l}(\bma,n)$ for the full range ${\bm k}, {\bm l}\in\Theta_n$, 
with arbitrary parameters $\bma$ and $n$. In Section~\ref{S:EE}, we describe such an algorithm 
that has the computational complexity order proportional to the number of these quantities. 

\section{Computing the B\'{e}zier coefficients $e^{\bm k}_{\bm l}(\bma,n)$}
					\label{S:EE}

The coefficients $e^{\sbm k}_{\sbm l}\equiv e^{\sbm k}_{\sbm l}(\bma,n)$, ${\bm k}, {\bm l}\in\Theta_n$, 
can be arranged in a block triangular
matrix  
\begin{equation}
	\label{E:E-tab}	
	\mathbb E=\left[\,\mathbb E^{\sbm k}\,\right]_{\sbm k\in\Theta_n}=\left[\begin{array}{lllll}
	\mathbb E^{(0,n)}\\[0.75ex]
	\mathbb E^{(0,n-1)}&\mathbb E^{(1,n-1)}\\[0.75ex]
\multicolumn{3}{c}{\dotfill}\\[0.75ex]
\mathbb E^{(0,1)}&\mathbb E^{(1,1)}&\ldots&\mathbb E^{(n-1,1)}\\[0.75ex]
\mathbb E^{(0,0)}&\mathbb E^{(1,0)}&\ldots&\mathbb E^{(n-1,0)}&\mathbb E^{(n,0)}
\end{array}\right],
\end{equation}
each block
$\mathbb E^{\sbm k}$, $\bm k\in\Theta_n$, being a triangular matrix 
\[  
	\mathbb E^{\sbm k}=\left[\,e^{\sbm k}_{\sbm l}\,\right]_{\sbm l\in\Theta_n}=
\left[\begin{array}{lllll}
e^{\sbm k}_{(0,n)}\\[1ex]
e^{\sbm k}_{(0,n-1)}&e^{\sbm k}_{(1,n-1)}\\[1ex]
\multicolumn{3}{c}{\dotfill}\\[1ex]
e^{\sbm k}_{(0,1)}&e^{\sbm k}_{(1,1)}&\ldots&e^{\sbm k}_{(n-1,1)}\\[1ex]
e^{\sbm k}_{(0,0)}&e^{\sbm k}_{(1,0)}&\ldots&e^{\sbm k}_{(n-1,0)}&e^{\sbm k}_{(n,0)}
\end{array}\right].
\]  	

The proposed algorithm is based on some recurrences satisfied by the elements of the matrix \eqref{E:E-tab}. In the sequel, we assume that $e^{\sbm k}_{\sbm l}=0$ if $\bm k\not\in\Theta_n$ 
or $\bm l\not\in\Theta_n$.

The \textit{first recurrence} is the following (cf.\ Lemma~\ref{L:e-recur}):
\begin{equation}
	\label{E:firstrec}
	e^{\sbm k+\sbm v_2}_{\sbm l}=\left([\sigma_1(\bm k)-\sigma_1(\bm l)]e^{\sbm k}_{\sbm l}
	-\sigma_2(\bm k) e^{\sbm k-\sbm v_2}_{\sbm l} 
	+\sigma_0(\bm l) e^{\sbm k}_{\sbm l+\sbm v_2}
	+\sigma_2(\bm l) e^{\sbm k}_{\sbm l-\sbm v_2}\right)/\sigma_0(\bm k),	
	\end{equation}
where $\bm v_2:=(0,1)$, and where for $\bm t:=(t_1,t_2)$, we define
\[ 	
	\sigma_0(\bm t):=(|\bm t|-n)(t_2+\alpha_2+1),\quad	
	\sigma_2(\bm t):=t_2(|\bm t|-\alpha_3-n-1),\quad
	\sigma_1(\bm t):=\sigma_0(\bm t)+\sigma_2(\bm t).
\]
Observe that recurrence \eqref{E:firstrec} relates three consecutive blocks of a column 
of the matrix \eqref{E:E-tab},  shown in the following diagram:
\[ 	
	\begin{array}{l}
	\mathbb E^{\sbm k+\sbm v_2}\\[0.75ex]
        \mathbb E^{\sbm k}\\[0.75ex]
	\mathbb E^{\sbm k-\sbm v_2}\,.
\end{array}
\] 
According to the  convention, the block $\mathbb E^{(k_1,-1)}$ has only zero elements.
Thus, we can compute all the blocks of this column, provided that we have computed the block 
$\mathbb E^{(k_1,0)}$ using another method.

Now, to \textit{initialise} the computation of the first column, we need a method to compute the corner block 
$\mathbb  E^{(0,0)}$ in the table~\eqref{E:E-tab}. We can use the following formula
(cf.\ Corollary~\ref{C:e00}):
\begin{equation}
	\label{E:E00}
	e^{(0,0)}_{\sbm l}(\bma,n) =\frac{(-1)^{l_1}(|\bma|+3)_{n}}{n!(\alpha_1+2)_{l_1}}
	\sum_{i=0}^{n-l_1}C^\ast_i\,h_i(l_2;\alpha_2,\alpha_3,n-l_1),
	\qquad \bm l:=(l_1,l_2)\in\Theta_n,
\end{equation}
where the coefficients $C^\ast_i$ are given by \eqref{E:Cast},  and we use the notation
$h_i(t;a,b,M)$  for the univariate Hahn polynomials (cf.~\eqref{E:Hahn1}).		
As noticed in Remark~\ref{R:Clenshaw}, we can efficiently  evaluate the sum
in \eqref{E:E00} using the Clenshaw's algorithm, at the cost of $O(n-l_1)$ operations.

To compute the remaining part of the last row of the matrix~\eqref{E:E-tab}, i.e.,
			\[
	\mathbb E^{(1,0)}\quad\mathbb E^{(2,0)}\quad\ldots\quad\mathbb E^{(n-1,0)}\quad\mathbb E^{(n,0)},
			\]	
we use the \textit{second recurrence} (cf.\ Lemma~\ref{L:e-recur}):
\begin{equation}\label{E:secrec}
	e^{\sbm k+\sbm v_1}_{\sbm l}
	=\left([\tau_1(\bm k)- \tau_1(\bm l)]e^{\sbm k}_{\sbm l}
	-\tau_2(\bm k)\,e^{\sbm k-\sbm v_1}_{\sbm l}
				 +
				 \tau_0(\bm l)\,e^{\sbm k}_{\sbm l+\sbm v_1}				  
				 + \tau_2(\bm l)\,e^{\sbm k}_{\sbm l-\sbm v_1}\right)/\tau_0(\bm k),		
\end{equation}
where 
$\bm v_1:=(1,0)$, 
and for   $\bm t:=(t_1,t_2)$, the coefficients  
$\tau_{j}(\bm t)$ are given by
\[  
	\tau_0(\bm t):=(|\bm t|-n)(t_1+\alpha_1+1),\quad	
	\tau_2(\bm t):=t_1(|\bm t|-\alpha_3-n-1), \quad 
	\tau_1(\bm t):=\tau_0(\bm t)+\tau_2(\bm t). 
\]  
The recurrence \eqref{E:secrec} relates each three consecutive blocks of  a row of the matrix \eqref{E:E-tab}, 
shown in the following diagram:
\[ 	
\begin{array}{lll}
\mathbb E^{\sbm k-\sbm v_1}&\quad\mathbb E^{\sbm k}&\quad\mathbb E^{\sbm k+\sbm v_1}\,.
\end{array}
\]  
Again, by the  convention, the blocks $\mathbb E^{(-1,k_2)}$ have only zero elements.

The following algorithm can be applied to compute the complete set of the coefficients
$e^{\sbm k}_{\sbm l}$ arranged in the block matrix $\mathbb E$ (cf.\ \eqref{E:E-tab}). The algorithm
requires $O(n^4)$ operations, i.e., its cost is proportional to the total
number of the coefficients.
			
\begin{alg}[Computing the table $\mathbb E$]
	                           \label{A:E-comp}
\
\begin{description}
\itemsep4pt  
\item[{\sc Step 1}]For $l_1=0,1,\ldots, n-1$,\\
\hphantom{Step 1 }$l_2=0,1,\ldots,n-l_1$,\\
\hphantom{Step 1 } using the Clenshaw algorithm, compute $e^{(0,0)}_{(l_1,l_2)}$ defined by \eqref{E:E00}, \eqref{E:Cast},\\
\hphantom{Step 1 } and put $e^{(l_1,l_2)}_{(0,0)}:=e^{(0,0)}_{(l_1,l_2)}$.

\item[{\sc Step 2}] For $k_1=0,1,\ldots,n-1$,
	\vspace*{-3ex}
	\begin{description}
\itemsep4pt    
\item[$1^o$]
	\hphantom{$1^o$} 
	for $k_2=0,1,\ldots,n-k_1-1$,\\
	\hphantom{$1^o$}$l_1=k_1,k_1+1,\ldots,n$,  \\
	\hphantom{$1^o$}$l_2=0,1,\ldots,n-l_1$,  \\
    \hphantom{$1^o$}
    compute $e^{(k_1,k_2+1)}_{(l_1,l_2)}$ using the recurrence \eqref{E:firstrec},
    and put $e^{(l_1,l_2)}_{(k_1,k_2+1)}:=e^{(k_1,k_2+1)}_{(l_1,l_2)}$;
\item[$2^o$]for $l_1=k_1+1,k_1+2,\ldots,n$,  \\
    \hphantom{$2^o$}$l_2=0,1,\ldots,n-l_1$,  \\
    \hphantom{$2^o$}
    compute $e^{(k_1+1,0)}_{(l_1,l_2)}$ using the recurrence \eqref{E:secrec},
    and put $e^{(l_1,l_2)}_{(k_1+1,0)}:=e^{(k_1+1,0)}_{(l_1,l_2)}$.
\end{description}
\end{description}	
\noindent \textbf{Output}: The set of the coefficients $e^{\sbm k}_{\sbm l}(\bma, n)$
for $\bm k,\bm l\in\Theta_n$.
\end{alg}

\begin{rem}\label{R:sym} 
	\begin{description}
	\itemsep1pt
	\item[(i)] In Algorithm~\ref{A:E-comp}, we made use of the symmetry property
	$e^{\sbm k}_{\sbm l}(\bma,n)=e^{\sbm l}_{\sbm k}(\bma,n)$
	(cf.\ \eqref{E:sym1}).
	\item[(ii)] When $\alpha_2=\alpha_3$, 
	the cost of completing the table $\mathbb E$ can be  reduced significantly
	using \eqref{E:sym2}. 
	First, the lower part of the table \eqref{E:E-tab}, containing the  blocks
        \[                                                             
	\left.\begin{array}{l}
	\mathbb E^{(k_1,\lfloor\frac12(n-k_1)\rfloor)} \\[1ex]
	\vdots\\[1ex]
	\mathbb E^{(k_1,1)}\\[1ex]
	\mathbb E^{(k_1,0)}
	\end{array}\;\right\}\;\mbox{with}\; k_1=0,1,\ldots,n,
	\]
	is computed using the properly adapted Algorithm~\ref{A:E-comp}.
	In the second step,  the upper part of the table is computed, using the formula 
	\(
		e^{\sbm k}_{\sbm l}(\bma,n)=e^{\hat{\sbm k}}_{\hat{\sbm l}}(\bma,n),
		\)
	where $\hat{\bm k}:=(k_1,n-|\bm k|)$ and $\hat{\bm l}:=(l_1,n-|\bm l|)$.
	\item[(iii)] Similar effect of the cost reduction can be obtained if 
			$\alpha_1=\alpha_2$, or
		$\alpha_1=\alpha_3$, or $\alpha_1=\alpha_2=\alpha_3$  
		(cf.\ \eqref{E:sym4}--\eqref{E:sym6}).
	\item[(iv)] Observe that the complexity order of Algorithm~\ref{A:E-comp} 
		equals   $O(n^4)$, i.e., is proportional to the total number of the coefficients 
		$e^{\sbm k}_{\sbm l}(\bma, c,n)$ for $\bm k,\bm l\in\Theta_n$.
\end{description}
\end{rem}
\section{Applications}
						\label{S:Appl}
Rational B\'ezier curves and surfaces, being the most natural generalization of polynomial
B\'ezier curves and surfaces, are important tools in geometric modelling.
However, they are sometimes inconvenient in practical applications. For this reason,
several algorithms for approximating a 
rational  B\'ezier geometric form by a polynomial one have been proposed 
(see, e.g., \cite{Hu13,LWK12,Rab05,Sha13}).

The following constrained approximation problem was recently considered in \cite{LKW17}.
\begin{prob}\label{P:approx}
Given the B\'{e}zier coefficients $r_{\sbm k}$ 
and positive weights $\omega_{\sbm k}$, $\bm k\in\Theta_n$, of the  rational  function $\sR_n$ 
of degree $n$,
\begin{equation}\label{E:Rn}
	\sR_n(\bm x)	
	:=\frac{\displaystyle \sum_{\sbm k\in\Theta_n}\omega_{\sbm k}r_{\sbm k}B^n_{\sbm k}(\bm x)}
	{\displaystyle \sum_{\sbm k\in\Theta_n}\omega_{\sbm k}B^n_{\sbm k}(\bm x)}
	= \sum_{\sbm k\in\Theta_n}r_{\sbm k}Q^n_{\sbm k}(\bm x),\qquad \bm x\in T,
\end{equation}
where
\begin{equation}
	\label{E:Q}
	Q^n_{\sbm k}(\bm x):=\frac{\displaystyle \omega_{\sbm k}B^n_{\sbm k}(\bm x)}
	{\displaystyle \sum_{\sbm i\in\Theta_n}\omega_{\sbm i}B^n_{\sbm i}(\bm x)},	
	\end{equation}
find a~polynomial of degree $m$, of the form
\begin{equation}		\label{E:Pm}
	\sP_m(\bm x)
	:= \sum_{\sbm k\in\Theta_m}p_{\sbm k}B^m_{\sbm k}(\bm x),\qquad \bm x\in T,	         
\end{equation}
with the coefficients $p_{\sbm k}$ 	satisfying the conditions
\begin{equation}\label{E:gcond}
	p_{\sbm k}=g_{\sbm k}\quad \mbox{for}\quad \bm k\in\Gamma^{\sbm c}_m,
\end{equation}	
$g_{\sbm k}$ being prescribed numbers, and  $\bm c:=(c_1,c_2,c_3)\in\N^3$ being a given
parameter vector with $|\bm c|<m$,
such that the distance between  $\sR_n$ and $\sP_m$,
\begin{equation}
	\label{E:dist}
	d(\sR_n,\sP_m):=\intT w_{\sbma}(\bm x)[\sR_n(\bm x)-\sP_m(\bm x)]^2\dx,
	\end{equation}
reaches the minimum.
\end{prob}

It has been shown that the solution of Problem~\ref{P:approx} is the polynomial \eqref{E:Pm}
with the coefficients given by (cf.\ \cite[Thm 2.2]{LKW17})
\begin{equation}\label{E:psc}
	p_{\sbm k}=\sum_{\sbm l\in\Omega^{\sbm c}_m}\binom{m}{\bm l}\,E^{\sbm k}_{\sbm l}(\bma,\bmc,m)
\big(u_{\sbm l}-v_{\sbm l}\big), \qquad \bm k\in\Omega^{\sbm c}_m,	
\end{equation}
where
\begin{align*}  
	u_{\sbm l}:=& \sum_{\sbm h\in\Theta_n} \binom{n}{\bm h}
	\rbinom{n+m}{\bm h+\bm l}\,r_{\sbm h}\,I_{\sbm h,\sbm l},\\ 
        \label{E:vl}
	v_{\sbm l}:=& \frac{1}{(|\bma|+3)_{2m}}\sum_{\sbm h\in\Gamma^{\sbm c}_m}\binom{m}{\bm h}	     	       
	       \left(\prod_{i=1}^{3}(\alpha_i+1)_{h_i+l_i}\right)g_{\sbm h}           	
\end{align*}
with $h_3:=m-|\bm h|$, $l_3:=m-|\bm l|$, and
\begin{equation}
	\label{E:I}
	I_{\sbm h,\sbm l}:=  \intT w_{\sbma}(\bm x)Q^{n}_{\sbm h}(\bm x)
	B^m_{\sbm l}(\bm x)\dx.
\end{equation}         
The symbol $E^{\sbm k}_{\sbm l}(\bma,\bmc,m)$ has the meaning given in \eqref{E:Dc2inB2-in}
(see also \eqref{E:Dc2inB2-coef-impl}).

Implementation of formula \eqref{E:psc} demands
$1^o$ evaluation of all  the coefficients $E^{\sbm k}_{\sbm l}(\bma,\bmc,m)$ 
with $\bm k, \bm l \in\Omega^{\sbmc}_m$ and $2^o$ computing all the integrals $I_{\sbm h,\sbm l}$
with $\bm h\in\Theta_n$, $\bm l\in\Omega^{\sbmc}_m$. 

The first task is accomplished in two steps. 
(\textit{i}) We compute all the the coefficients $e^{\sbm k}_{\sbm l}(\bm\mu,\bmc,M)$,
$\bm k, \bm l \in\Omega^{\sbmc}_m$, with $M:=m-\lbmc$ and $\bm\mu:=\bma+2\bmc$
by Algorithm~\ref{A:E-comp}
and then (\textit{ii}) use the result in formula \eqref{E:Dc2inB2-coef-impl} (with $n$
replaced by $m$) to compute 
$E^{\sbm k}_{\sbm l}(\bma,\bmc,m)$, $\bm k, \bm l \in\Omega^{\sbmc}_m$. 
Remark that the striking simplicity of Algorithm~\ref{A:E-comp} was obtained  
thanks to using  properties  of dual bivariate Bernstein polynomials, investigated in Appendix~A.

As for the second task, observe that in general,
the integrals \eqref{E:I} cannot be evaluated exactly,  and their
	number, equal to $(n+m-|\bm c|)(n+m-|\bm c|+1)/2$, may lead to the impression 
	that computation of the coefficients \eqref{E:psc} is time-expensive in practice. 
    The problem has been defeated in \cite{LKW17},
    where we have proposed a very effective algorithm which allows to compute the complete set
    of integrals $I_{\sbm h,\sbm l}$, $\bm h\in\Theta_n$, $\bm l\in\Omega^{\sbm c}_m$, in a time
    about twice as long as the time required to compute a single integral of this type (the
    algorithm is an extension of the adaptive, high precision method of \cite{Kel07},
    and it uses some special properties of the integrals \eqref{E:I}).
       
    In this way, we have obtained a very efficient method for computing the solution to the
Problem~\ref{P:approx}. 

    \subsection{ Case of equal weights}
    \label{SS:eqw}

Notice that in the particular case where all the weights are equal, $\omega_{\sbm i}=\omega$, $\bm i\in\Theta_n$,
we have $Q^n_{\sbm k}(\bm x)=B^n_{\sbm k}(\bm x)$, $\bm k\in\Theta_n$. Consequently, the rational function
	\eqref{E:Rn} reduces to a polynomial of degree $n$, so that  
	Problem~\ref{P:approx} (with additional assumption $m<n$) 
	is actually the constrained polynomial degree reduction problem 
	that has been discussed \cite{WL10}. 
	In the cited paper, we have developed a method to evaluate the control points
of the reduced surface in terms of quantities 
	$\left\langle D^{(m,\sbm c)}_{\sbm k},\,B^n_{\sbm l} \right\rangle_{\sbm a}$, 
	$\bm k\in \Omega^{\sbm c}_{m}$, $\bm l\in\Theta_n$.
	Computing these coefficients was the most cumbersome part of the method. 
The main part of the above collection
was computed efficiently thanks to using some recurrence relations. However, 
a quite large portion of the coefficients  demanded some complex computations. 
By applying the approach of this paper, we may solve the degree reduction problem much faster.

N.\,B. In the case of equal weights, evaluation of the  integrals \eqref{E:I} is a simple task as we have
(see, e.g., \cite{WL10})
	\[
		I_{\sbm h,\sbm l}=\binom{n}{\bm h}\binom{m}{\bm l}
	\frac{
	(\alpha_1+1)_{h_1+l_1}(\alpha_2+1)_{h_2+l_2}(\alpha_3+1)_{n+m-|\sbm h|-|\sbm l|}}
	{(|\bma|+3)_{n+m}}.	
		\]

\section*{Appendix A: Dual bivariate Bernstein bases}
				\label{S:AppA}
\setcounter{equation}0
\setcounter{subsection}0
\setcounter{thm}0
\renewcommand{\thethm}{A.\arabic{thm}}
\renewcommand{\thesubsection}{A.\arabic{subsection}}
\renewcommand{\theequation}{A.\arabic{equation}}

\subsection{Unconstrained dual bivariate Bernstein polynomials}
\label{SS:dualB2}

In this section, we prove some important properties of the coefficients $e^{\sbm k}_{\sbm l}(\bma,n)$
introduced in \eqref{E:D2inB2-in}. Obviously, we have (cf.\ \eqref{E:Bezform})
\begin{equation}\label{E:e-is}
	e^{\sbm k}_{\sbm l}(\bma,n)=\langle D^n_{\sbm k}, D^n_{\sbm l}\rangle_{\sbma}, \qquad
	\bm k,\bm l\in\Theta_n.	
	\end{equation}
However, we need some alternative formulas that would be more  useful in obtaining recurrence relations 
satisfied by these quantities.

In the sequel, we adopt the following convention: given $\bm t:=(t_1,t_2)\in\Theta_n$, we use the notation 
\[
	\tilde{\bm t}:=(t_2,t_1),\quad \hat{\bm t}:=(t_1,t_3),\quad \check{\bm t}:=(t_3,t_1),\quad \bm t^\ast:=(t_2,t_3)\quad \mbox{and}\quad {\bm t}^{\circ}:=(t_3,t_2),\
	\]
	where $t_3:=n-|\bm t|$.	
 \begin{lem}\label{L:D2inB2}
 	 Dual bivariate Bernstein polynomials have the following representation:
\[ 
	D^n_{\sbm k}(\bm x;\bma)=\,
	\sum_{\sbm l\in\Theta_n}\,e^{\sbm k}_{\sbm l}(\bma,n)\,
	B^n_{\sbm l}(\bm x),\qquad \bm k\in\Theta_n,
\]  
where
\begin{equation}\label{E:D2inB2-coef}
	e^{\sbm k}_{\sbm l}(\bma,n):=\sum_{0\le i\le q\le n}\,C^2_{q,i}(\bma,n)
	H_{q,i}\left(\bm k^\ast;\bma,n\right)H_{q,i}\left(\bm l^\ast;\bma,n\right)
\end{equation}
with 
\[  
	C_{q,i}(\bma,n):=\binom{q}{i}\frac{1}{q!(-n)_q\, \lambda_{q,i}}.
\]
Here $H_{q,i}(\bm t;\bma,n)$ are the bivariate Hahn polynomials defined by \eqref{E:Hahn2},
and $\lambda_{q,i}$ is the constant given by \eqref{E:lam}.
\end{lem}
\begin{pf}
	In \cite{LW06}, it has been proved that the dual bivariate Bernstein polynomials 
have the following  representation:
\begin{equation}\label{E:D2inP2}
	D^n_{\sbm k}(\bm x;\bma)=\,
	\sum_{0\le i\le q\le n}\,b^{\sbma}_{q,i}(n,\bm k)\,
	P_{q,i}^{\sbma}(\bm x),\qquad \bm k\in\Theta_n,
\end{equation}
where
\[  
	b^{\sbma}_{ q,i}(n,\bm k):=
	(-1)^{i}C_{q,i}(\bma,n)
	H_{q,i}\left(\bm k^\ast;\bma,n\right).
\] 
Here  $P_{q,i}^{\sbma}$,  $q\in\N;\; 0\le i\le q$,
are the \textit{two-variable Jacobi polynomials}
(\!\!\cite{Koo75}; see also \cite[p. 86]{DX01} or \eqref{E:Jacxy}), which 
 are orthonormal on $T$, i.e.,  
 $\langle P_{m,l}^{\sbma},\, P_{n,k}^{\sbma} \rangle_{\sbma}$ equals
 1 if $(m,l)=(n,k)$,  and 0 otherwise.
Using  expansion \eqref{E:D2inP2} for both dual polynomials in \eqref{E:e-is}, we obtain
\[
	e^{\sbm k}_{\sbm l}(\bma,n)
	=\sum_{0\le i\le q\le n}b^{\sbma}_{q,i}(n,\bm k)b^{\sbma}_{q,i}(n,\bm l),
	\]
	which is equivalent to \eqref{E:D2inB2-coef}.	
	\end{pf}

\begin{rem} \label{R:ekl-sym}
	Some symmetry properties of the coefficients $e^{\sbm k}_{\sbm l}(\bma,n)$  should  be noticed. 
\begin{enumerate}
\itemsep 2pt		
\item   The following equation results easily from \eqref{E:e-is}:	
\begin{equation}
	\label{E:sym1}
	e^{\sbm k}_{\sbm l}(\bma,n)=e^{\sbm l}_{\sbm k}(\bma,n).		
	\end{equation}	

\item  Let $\alpha_2=\alpha_3$. The following equality holds:
	\begin{equation}
		\label{E:sym2}
		e^{\sbm k}_{\sbm l}(\bma,n)=e^{\hat{\sbm k}}_{\hat{\sbm l}}(\bma,n).
		\end{equation}
	This equation  can be easily verified using 	\eqref{E:D2inB2-coef},  definition \eqref{E:Hahn2}
	of the bivariate Hahn polynomials, and the identity 
		\(
			h_i(n-|\bm t|;\alpha_2,\alpha_2,n-t_1)=(-1)^ih_i(t_2;\alpha_2,\alpha_2,n-t_1)
			\) (cf. \eqref{E:Hahn1}).
\item 	The following equation holds:
	\begin{equation}
	\label{E:sym3}
	e^{\sbm k}_{\sbm l}(\bma,n)=e^{\tsbm k}_{\tsbm l}(\tilde{\bma},n),
	\end{equation}
where  $\tilde{\bma}:=(\alpha_2,\alpha_1,\alpha_3)$.
By \eqref{E:e-is}, equation \eqref{E:sym3} is equivalent to the  equation 
\[
	\left\langle   D^n_{\sbm k}(\cdot;\bma),D^n_{\sbm l}(\cdot;\bma) \right \rangle_{\sbma}
	= \left\langle D^n_{\tsbm k}(\cdot;\tilde{\bma}),D^n_{\tsbm l}(\cdot;\tilde{\bma})\right\rangle_{\tsbm{\alpha}}
	\]
	which can be easily verified using the definition \eqref{E:Jinprod} of the 
	inner product $\langle\cdot,\cdot\rangle_{\sbma}$.				
\item  	Let $\alpha_1=\alpha_2$. By \eqref{E:sym3}, we obtain
			the equation
	\begin{equation}
		\label{E:sym4}
		e^{\sbm k}_{\sbm l}(\bma,n)=e^{\tsbm k}_{\tsbm l}(\bma,n).
\end{equation}

\item  	Let $\alpha_1=\alpha_3$. Then
\begin{equation}
		\label{E:sym5}
		e^{\sbm k}_{\sbm l}(\bma,n)=e^{{\sbm k}^{\circ}}_{{\sbm l}^{\circ}}(\bma,n).
\end{equation}
Transforming $e^{\sbm k}_{\sbm l}(\bma,n)$, using consecutively \eqref{E:sym3}, \eqref{E:sym2}, and again \eqref{E:sym3}, the result follows. 
\item Let $\alpha_1=\alpha_2=\alpha_3$. By \eqref{E:sym2}, \eqref{E:sym4}, and \eqref{E:sym5},
\begin{equation}
		\label{E:sym6}
		e^{\sbm k}_{\sbm l}(\bma,n)=e^{\tsbm k}_{\tsbm l}(\bma,n)=e^{\hat{\sbm k}}_{\hat{\sbm l}}(\bma,n)
		=e^{\sbm k^\ast}_{\sbm l^\ast}(\bma,n)=e^{\check{\sbm k}}_{\check{\sbm l}}(\bma,n)
		=e^{{\sbm k}^{\circ}}_{{\sbm l}^{\circ}}(\bma,n).		
\end{equation}
\end{enumerate}
\end{rem}

\begin{lem}\label{L:e-recur}
	The coefficients  $e^{\sbm k}_{\sbm l}\equiv e^{\sbm k}_{\sbm l}(\bma,n)$
    defined in \eqref{E:D2inB2-in}, and
    given by \eqref{E:D2inB2-coef},
	satisfy the following bivariate recurrence relations:
\begin{align}
	\label{E:e-receq1}
			&\sigma_0(\bm k)\, e^{\sbm k+\sbm v_2}_{\sbm l}-\sigma_1(\bm k)\, e^{\sbm k}_{\sbm l}
			+\sigma_2(\bm k) \,e^{\sbm k-\sbm v_2}_{\sbm l}			
			=\sigma_0(\bm l)\, e^{\sbm k}_{\sbm l+\sbm v_2}-\sigma_1(\bm l) \,e^{\sbm k}_{\sbm l}
			+\sigma_2(\bm l)\, e^{\sbm k}_{\sbm l-\sbm v_2},\\[1ex]	
	\label{E:e-receq2}
	& \tau_0(\bm k)\,e^{\sbm k+\sbm v_1}_{\sbm l}
				  -\tau_1(\bm k)\,e^{\sbm k}_{\sbm l}
				  +\tau_2(\bm k)\,e^{\sbm k-\sbm v_1}_{\sbm l}  
				=
				  \tau_0(\bm l)\,e^{\sbm k}_{\sbm l+\sbm v_1}
				  - \tau_1(\bm l)\,e^{\sbm k}_{\sbm l}
				  + \tau_2(\bm l)\,e^{\sbm k}_{\sbm l-\sbm v_1},	  
				  \end{align}			 
 where  $\bm v_1:=(1,0)$,  $\bm v_2:=(0,1)$, and for $\bm t:=(t_1,t_2)$, we define 
\begin{equation}
	\label{E:sigtau}
	\sigma_i(\bm t):=\varphi_i(\bma,\bm t),\quad \tau_i(\bm t):=\varphi_i(\tilde\bma,\tilde{\bm t}),\qquad i=0,1,2,	
	\end{equation}
with $\tilde\bma:=(\alpha_2,\alpha_1,\alpha_3)$,   and
\begin{equation}
	\label{E:varphi}
	\begin{array}{c}
	\varphi_0(\bma,\bm t):=(|\bm t|-n)(t_2+\alpha_2+1),\quad	
	\varphi_2(\bma,\bm t):=t_2(|\bm t|-\alpha_3-n-1),\\[1ex]
	\varphi_1(\bma,\bm t):=\varphi_0(\bma,\bm t)+\varphi_2(\bma,\bm t).
	\end{array}
	\end{equation}	
\end{lem}	
\begin{pf} First, we prove the recurrence \eqref{E:e-receq1}. By \eqref{E:Hahn2}, we have
	\[
	H_{q,i}\left(\bm t^\ast;\bma,n\right)=h_i(t_2;\alpha_2,\alpha_3,n-t_1)\,
	h_{q-i}(n-t_1-i;\alpha_2+\alpha_3+2i+1,\alpha_1,n-i).
		\]
		Let $\sD^{n-t_1}_{t_2}$ be the difference operator defined according to \eqref{E:Dop}.
		Then, by \eqref{E:Hahn1deq},
		\[
			\sD^{n-t_1}_{t_2}\,h_i(t_2;\alpha_2,\alpha_3,n-t_1)
			=i(i+\alpha_2+\alpha_3+1)\,h_i(t_2;\alpha_2,\alpha_3,n-t_1).	
			\]	
			Hence, we obtain
	\begin{align*}
				\sD^{n-t_1}_{t_2}\,H_{q,i}\left(\bm t^\ast;\bma,n\right)&=
			\sD^{n-t_1}_{t_2}\,h_i(t_2;\alpha_2,\alpha_3,n-t_1)\cdot 
			h_{q-i}(n-t_1-i;\alpha_2+\alpha_3+2i+1,\alpha_1,n-i)\\
			&=i(i+\alpha_2+\alpha_3+1) \,H_{q,i}\left(\bm t^\ast;\bma,n\right).
	  \end{align*}
	Consequently, having in mind the form \eqref{E:D2inB2-coef}, we obtain the equation
        \[
        	\sD^{n-k_1}_{k_2}e^{\sbm k}_{\sbm l}(\bma,n) = \sD^{n-l_1}_{l_2}e^{\sbm k}_{\sbm l}(\bma,n)
        \]	
		which can be simplified to the form \eqref{E:e-receq1}.

		To prove \eqref{E:e-receq2}, let us define the difference operator $\cR^{\sbma}_{\sbm t}$ by
	\[
		\cR^{\sbma}_{\sbm t}\,F(\bm t)=\varphi_0(\bma,\bm t)F(\bm t+\bm v_2)
		-\varphi_1(\bma,\bm t)F(\bm t)
		+\varphi_2(\bma,\bm t)F(\bm t-\bm v_2),
		\]
		where we use the notation of \eqref{E:varphi}.		
		The recurrence \eqref{E:e-receq1} can be written as 
	\[
		\cR^{\sbma}_{\sbm k}e^{\sbm k}_{\sbm l}(\bma,n)=
		\cR^{\sbma}_{\sbm l}e^{\sbm k}_{\sbm l}(\bma,n).			
		\]
		Substituting $\tilde{\bm k}$, $\tilde{\bm l}$ and $\tbma$ in place of $\bm k$, $\bm l$ and $\bma$,
		respectively, gives	
		\begin{equation}\label{E:Rt}
				\cR^{\tsbma}_{\tsbm k}e^{\tsbm k}_{\tsbm l}(\tbma,n)=
			\cR^{\tsbma}_{\tsbm l}e^{\tsbm k}_{\tsbm l}(\tbma,n).						
			\end{equation}
		Now, notice that by \eqref{E:sym3} we have
		\begin{align*}
				  & 
				  \cR^{\tsbma}_{\tsbm k}e^{\tsbm k}_{\tsbm l}(\tilde{\bma},n)
				  =\varphi_0(\tilde{\bma},\tilde{\bm k})e^{\sbm k+\sbm v_1}_{\sbm l}
				  -\varphi_1(\tilde\bma,\tilde{\bm k})e^{\sbm k}_{\sbm l}
				  +\varphi_2(\tilde\bma,\tilde{\bm k})e^{\sbm k-\sbm v_1}_{\sbm l}, \\
				  &\cR^{\tsbma}_{\tsbm l}e^{\tsbm k}_{\tsbm l}(\tilde{\bma},n)
				  =\varphi_0(\tilde\bma,\tilde{\bm l})e^{\sbm k}_{\sbm l+\sbm v_1}
				  -\varphi_1(\tilde\bma,\tilde{\bm l})e^{\sbm k}_{\sbm l}
				  +\varphi_2(\tilde\bma,\tilde{\bm l})e^{\sbm k}_{\sbm l-\sbm v_1}.
			\end{align*}
		Using this in \eqref{E:Rt}, equation \eqref{E:e-receq2} follows.		
	\end{pf}

\begin{lem}\label{L:e-k20}
	For $\bm k=(k_1,0),\,\bm l=(l_1,l_2)\in\Theta_n$, the following formula holds:
	\begin{equation}\label{E:e-k20}
	e^{\sbm k}_{\sbm l}(\bma,n) =
	 G_{\sbma}\sum_{i=0}^{n-\max(k_1,l_1)}C_i\,h_i(l_2;\alpha_2,\alpha_3,n-l_1),
			\end{equation}
		where $ G_{\sbma} :={\Gamma(\eta+1)\Gamma(\alpha_1+1)}/{\Gamma(|\bma|+3)}$, $\eta:=\alpha_2+\alpha_3+1$, symbol $h_i(t;a,b,M)$ is defined in \eqref{E:Hahn1}, and
			\begin{equation}
				\label{E:Ci}
				\left.\begin{array}{ll}
				C_0:=c_{k_1,l_1}(n,\eta,\alpha_1),\\[1ex]
				C_i:=\dfrac{(2i+\eta)(k_1-n)_i(\eta+1)_{i-1}}{i!(-n)_i^2(\alpha_3+1)_i}
				\,c_{k_1,l_1}(n-i,2i+\eta,\alpha_1),&\quad i\ge1,	
				\end{array}\quad\right\}
				\end{equation}		  				
			the  symbol $c_{h,j}(m,\alpha,\beta)$ denoting the $j$\emph{th} 
			B\'ezier coefficient of the univariate dual Bernstein polynomial 
			$D^{m}_{h}(x;\alpha,\beta)$  (cf.\ \cite{LW11b}).		
\end{lem}
\begin{pf} We give a sketch of the proof.
	By \cite[Thm 3.3]{LW06}, we have for $\bm k:=(k_1,k_2)\in\Theta_n$,
	\[
	D^n_{\sbm k}(\bm x;\bma)
	=G_{\sbma}\sum_{i=0}^{n-k_1}f_i(n,\bm k)
	(1-x_1)^iR^{(\alpha_3,\alpha_2)}_i\left(x_2/(1-x_1)\right)D^{n-i}_{k_1}(x_1;2i+\eta,\alpha_1),
	\]
	where    $R^{(\alpha_3,\alpha_2)}_i$  are the univariate Jacobi polynomials (cf.~\eqref{E:uniJac}),
	$D^N_j(t;\mu,\nu)$ are univariate dual Bernstein polynomials (see, e.g., \cite{LW11b}), and
	\[
		\begin{array}{ll}
		f_0(n,\bm k):=1,\\[0.5ex]
		f_i(n,\bm k):=(-1)^i\dfrac{(2i+\eta)(\eta+1)_{i-1}}
				           {(-n)_i(\alpha_2+1)_i(\alpha_3+1)_i}\,
				           \,h_i(k_2;\alpha_2,\alpha_3,n-k_1),&\quad i\ge1.
				           \end{array}
	\]
	We use the above formula in \eqref{E:e-is},
	then reduce the integration over the triangle $T$ to evaluating two one-dimensional integrals,
	and apply the orthogonality property of the  polynomials
	$R^{(\alpha_3,\alpha_2)}_i$ (cf.\ \cite[\S1.8]{KS98}). 
	Letting $\bm k:=(k_1,0)$, the formula \eqref{E:e-k20} follows.
\end{pf}

\begin{cor}\label{C:e00}
	The following formula holds:
	\[
			e^{(0,0)}_{\sbm l}(\bma,n) =\frac{(-1)^{l_1}(|\bma|+3)_{n}}{n!(\alpha_1+2)_{l_1}}
	\sum_{i=0}^{n-l_1}C^\ast_i\,h_i(l_2;\alpha_2,\alpha_3,n-l_1),
	\qquad \bm l:=(l_1,l_2)\in\Theta_n,
	\]
		where
	\begin{equation}
				\label{E:Cast}
				\left.\begin{array}{ll}
				C^\ast_0:=\dfrac{(\alpha_1+2)_{n}}
				{(\alpha_2+\alpha_3+2)_{n-l_1}},\\[2.5ex]
				\displaystyle C^\ast_i:=(-1)^{i}\dfrac{(2i+\alpha_2+\alpha_3+1)(\alpha_1+2)_{n-i}(|\bma|+n+3)_{i}}
				{i!(\alpha_3+1)_i(\alpha_2+\alpha_3+i+1)_{n-l_1+1}},&\quad i\ge1.	
				\end{array}\quad\right\}
\end{equation}				
\end{cor}	
\begin{pf}
	The result follows by putting $k_1=0$ in \eqref{E:e-k20}, \eqref{E:Ci}, and 
	using  the explicit form for  $c_{0,l_1}(n-i,2i+\alpha_2+\alpha_3+1,\alpha_1)$, given in 
	\cite[Eq. (2.11)]{LW11b}.
\end{pf}

\subsection{Constrained dual bivariate Bernstein polynomials}
									\label{SS:cdualA2}

The  constrained  dual bivariate Bernstein basis polynomials \eqref{E:constrdBer2} can be
expressed in terms
of the unconstrained dual bivariate Bernstein polynomials  \eqref{E:dBer2}
with shifted  degree and parameters. Namely, we have the following result.
\begin{lem}[\!\!\cite{WL10}]
	  \label{L:constrDual2}
	For $\bm k\in \Omega^{\sbm c}_n$, the following formula holds:
\begin{equation}
	\label{E:Dc2inD2}
	D^{(n,\sbm c)}_{\sbm k}(\bm x;\bma)=
	U\,V_{\sbm k}(n)\,
	x_1^{c_1}x_2^{c_2}(1-|\bm x|)^{c_3}\,
	D^{n-|\sbm c|}_{\sbm k-\sbm c'}(\bm x;\bma+2\bm c),
\end{equation}
where $\bm c':=(c_1,c_2)$, and
\begin{equation}
	\label{E:UV}
	 U:=
	(|\bma|+3)_{2|\sbm c|}\prod_{i=1}^{3}
	(\alpha_i+1)_{2c_i}^{-1},\qquad
	V_{\sbm k}(n):= \binom{n-|\bm c|}{\bm k-\bm c'}\rbinom{n}{\bm k}.
  \end{equation}
\end{lem}

\begin{lem}\label{L:Dc2inB2}
	The constrained dual bivariate Bernstein polynomials have the B\'ezier  representation
	\[ 
	D^{(n,\sbm c)}_{\sbm k}(\bm x;\bma)=\sum_{\sbm l\in\Omega^{\sbm c}_n}
	\,E^{\sbm k}_{\sbm l}(\bma,\bmc,n)\,B^n_{\sbm l}(\bm x),
	\] 	
where
\[ 
	E^{\sbm k}_{\sbm l}(\bma,\bmc,n):=U\,V_{\sbm k}(n)\,V_{\sbm l}(n)\,
	\,e^{\sbm k-\sbm c'}_{\sbm l-\sbm c'}(\bma+2\bm c,n-|\bm c|),\qquad \bm c':=(c_1,c_2),	
\] 
notation used being that of \eqref{E:D2inB2-coef} and \eqref{E:UV}.
\end{lem}
\begin{pf}
	By Lemma~\ref{L:D2inB2}, we have	
	\[
		D^{n-|\sbm c|}_{\sbm k-\sbm c'}(\bm x;\bma+2\bm c)
		=\sum_{\sbm l\in\Theta_{n-|\tbm c|}}	
		e^{\sbm k-\sbm c'}_{\sbm l}(\bma+2\bm c,n-|\bm c|)\,
		B^{n-|\sbm c|}_{\sbm l}(\bm x).
	\]
	Putting this result in \eqref{E:Dc2inD2} and using  the equation 
	\[
		x_1^{c_1}x_2^{c_2}(1-|\bm x|)^{c_3}\cdot B^{n-|\sbm c|}_{\sbm l}(\bm x)=
		\binom{n-|\bm c|}{\bm l}\binom{n}{\bm l+\bm c'}^{\!-1}\,B^n_{\sbm l+\sbm c'}(\bm x),
	\]
	we obtain
\begin{align*}
	D^{(n,\sbm c)}_{\sbm k}(\bm x;\bma)=&
	U\,V_{\sbm k}(n)\,\sum_{\sbm l\in\Theta_{n}}	
	\binom{n-|\bm c|}{\bm l}\rbinom{n}{\bm l+\bm c'}
	\,e^{\sbm k-\sbm c'}_{\sbm l}(\bma+2\bm c,n-|\bm c|)\,
	B^{n}_{\sbm l+\sbm c'}(\bm x)\\
	=&
	U\,V_{\sbm k}(n)\,\sum_{\sbm l\in\Omega^{\sbm c}_{n}}
	\binom{n-|\bm c|}{\bm l-\bm c'}\rbinom{n}{\bm l}	
	\,e^{\sbm k-\sbm c'}_{\sbm l-\sbm c'}(\bma+2\bm c,n-|\bm c|)\,
	B^{n}_{\sbm l}(\bm x).
\end{align*}
Hence, the lemma  is proved.                                       
\end{pf}

\section*{Appendix B: Bivariate orthogonal polynomials}
			\label{S:AppB}

\setcounter{equation}0
\setcounter{subsection}0
\setcounter{thm}0
\renewcommand{\thethm}{B.\arabic{thm}}
\renewcommand{\thesection}{\Arabic{section}}
\renewcommand{\thesubsection}{B.\arabic{subsection}}
\renewcommand{\theequation}{B.\arabic{equation}}

The notation
\[
 \hyper rs {a_1,\ldots, a_r}{b_1,\ldots,b_s}z
 		:=\sum_{k=0}^{\infty}\frac{(a_1)_k\cdots (a_r)_k}
                           {k!(b_1)_k\cdots (b_s)_k}\,z^k
\]
is used for the \textit{generalised hypergeometric series}
(see, e.g., \cite[\S2.1]{AAR99});
here $r,\,s\in\Zplus$, $z$, $a_1,\ldots, a_r$, $b_1,\ldots,\,b_s$$\in \C$,
and $(c)_k$ is
the shifted factorial.

Recall that \textit{two-variable Jacobi polynomials} 
$P_{n,k}^{\sbma}(\bm x)$, $n=0,1,\ldots,\; k=0,1,\ldots ,n$, are defined by
(\!\!\cite{Koo75}; see also \cite[p. 86]{DX01})
\begin{align}\label{E:Jacxy}
	P_{n,k}^{\sbma}(\bm x):=& \lambda^{-1}_{n,k}
         \,R_{n-k}^{(2k+\alpha_2+\alpha_3+1,\alpha_1)}(x_1)
         \,(1-x_1)^k\,
         R_k^{(\alpha_3 ,\alpha_2)}\left({x_2}/{(1-x_1)}\right),
\end{align}
where $\bm x:=(x_1,x_2)$, $\bma:=(\alpha_1,\,\alpha_2,\,\alpha_3)$, $\alpha_i>-1$, 
\begin{equation}\label{E:uniJac}
R^{(\mu,\nu)}_m(t):=\frac{(\mu+1)_m}{m!}\hyper21{-m,m+\mu+\nu+1}{\mu+1}{1-t}
\end{equation}
is the $m$th shifted Jacobi polynomial in one variable \cite[\S 1.8]{KS98},
and
\begin{equation}\label{E:lam}
 \lambda^{2}_{n,k}\equiv
 \left[ \lambda^{\sbma}_{n,k}\right]^{2}:=
 \frac{(\alpha_1+1)_{n-k}
 (\alpha_2+1)_k(\alpha_3+1 )_k(k+\eta)_{n+1}}
 {k!(n-k)!(2k+\eta)(2n/\sigma+1)(\sigma)_{n+k}}
\end{equation}
with
 $\eta:=\alpha_2+\alpha_3+1$,  $\sigma:=|\bma|+2$.
 Polynomials \eqref{E:Jacxy} form  the orthonormal set with respect 
 to the inner product 
 \[  
	\langle f, g \rangle_{\sbma} :=\intT w_{\sbma}(\bm x)f(\bm x)\,g(\bm x) \dx,
 \] 
where  $T:=\{(x_1,\,x_2)\;:\;x_1,\,x_2\ge0,\; 1-x_1-x_2\ge0\}$, and 
$w_{\sbma}(\bm x):=A_{\sbma}\, x_1^{\alpha_1} x_2^{\alpha_2}(1-x_1-x_2)^{\alpha_3}$ with
\begin{equation}\label{E:A}
A_{\sbma}:=\Gamma(\lbma+3)\prod_{i=1}^{3}[\Gamma(\alpha_i+1)]^{-1}.
\end{equation}
\textit{Bivariate Hahn polynomials} are defined by (see, e.g., \cite{Tra91}) 
\begin{equation}\label{E:Hahn2}
	H_{q,i}(\bm t;\bma,n):=h_i(t_1;\alpha_2,\alpha_3,t_1+t_2)
	\,h_{q-i}(t_1+t_2-i;\alpha_2+\alpha_3+2i+1,\alpha_1,n-i)\,,
\end{equation}
where $0\le i\le q\le n$, $n\in\N$, $\bma=:(\alpha_1,\alpha_2\,\alpha_3)$, 
$\alpha_i>-1$, $i=1,2,3$, $\bm t:=(t_1,t_2)\in\R^2$, and
\begin{equation}
	\label{E:Hahn1}
	h_l(t;a,b,M):=(a+1)_l(-M)_l\hyper32{-l,l+a+b+1,-t}{a+1,-M}1
\end{equation}  
are the \textit{univariate Hahn polynomials} (see, e.g.,  \cite[\S1.5]{KS98}). The latter polynomials satisfy the  recurrence relation
\[ 
		h_{l+1}(t)=A_l(t,M)\,h_l(t)+B_l(M)\,h_{l-1}(t),
		\qquad l\ge0;\;	h_{0}(t)\equiv1;	\;h_{-1}(t)\equiv0,
\]  
 where $h_l(t)\equiv h_l(t;a,b,M)$,
\begin{equation}	\label{E:Hahn1-rec-coeffs}
	A_l(t,M):=C_l\,(2l+s-1)_2\,t-D_l-E_l,\qquad 	
 B_l(M):=-D_l\,E_{l-1},		  	
\end{equation}
with $s:=a+b+1$, $C_l:=(2l+s+1)/[(l+s)(2l+s-1)]$, $D_l:=C_l\,l(l+M+s)(l+b)$, and $E_l:=(l+a+1)(M-l)$.
Moreover,  the polynomials \eqref{E:Hahn1}
satisfy the difference equation
\begin{equation}\label{E:Hahn1deq}
	\sD^M_t\,h_j(t;a,b,M)=j(j+a+b+1)\,h_j(t;a,b,M),
	\end{equation}
where the difference operator $\sD^M_t$ is given by
\begin{equation}
	\label{E:Dop}
	\sD^M_t F(t):=U(t;a,b,M)\,F(t+1)-V(t;a,b,M)\,F(t)+W(t;a,b,M)\,F(t-1)
	\end{equation}
with $U(t;a,b,M):=(t-M)(t+a+1)$, $W(t;a,b,M):=t(t-b-M-1)$ and
	$V(t;a,b,M):=U(t;a,b,M)+W(t;a,b,M)$.

\begin{rem} \label{R:Clenshaw}
	A linear combination 
	\(
		s_N(t):=\sum_{i=0}^{N}\gamma_ih_i(t;a,b,M)
	\)	
	can be summed up using the following \textit{Clenshaw's algorithm} (see, e.g., \cite[Thm 3.2.11]{DB08}).				
	Compute the sequence $V_0,V_1,\ldots$, $V_{n+2}$ from				
	\[
		V_i:=\gamma_i+A_i(t;M)V_{i+1}+B_{i+1}(M)V_{i+2},\qquad i=N,N-1,\ldots,0,				
	\]
	with $V_{N+1}=V_{N+2}=0$, where the coefficients $A_i(t;M)$ and $B_i(M)$ are defined by 
	\eqref{E:Hahn1-rec-coeffs}.
	Then, $s_N(t)=V_0$.			 				 
\end{rem}


\end{document}